\documentclass[12pt,a4paper]{amsart}

\usepackage{amscd, amssymb}

\begin{document}

\pagestyle{headings}

\def\clsp{\overline{\operatorname{span}}}
\def\newspan{{\operatorname{span}}}

\newtheorem{thm}{Theorem}[section]
\newtheorem{cor}[thm]{Corollary}
\newtheorem{lem}[thm]{Lemma}
\newtheorem{prop}[thm]{Proposition}
\newtheorem{thm1}{Theorem}

\theoremstyle{definition}
\newtheorem{dfn}[thm]{Definition}
\newtheorem{dfns}[thm]{Definitions}

\theoremstyle{remark}
\newtheorem{rmk}[thm]{Remark}
\newtheorem{rmks}[thm]{Remarks}
\newtheorem{example}[thm]{Example}
\newtheorem{examples}[thm]{Examples}
\newtheorem{note}[thm]{Note}
\newtheorem{notes}[thm]{Notes}

\def\Ind{\operatorname{Ind}}
\def\rt{\operatorname{rt}}
\def\id{\operatorname{id}}
\def\supp{\operatorname{supp}}
\def\Aut{\operatorname{Aut}}
\numberwithin{equation}{section}
\def\H{\mathcal{H}}
\def\K{\mathcal{K}}

\title[Coverings of directed graphs]{\boldmath Coverings of directed 
graphs\\ and
crossed products of $C^*$-algebras\\ by coactions of homogeneous spaces}
\author{Klaus Deicke}
\address{Department of Mathematics, University of Paderborn, 33095 
Paderborn, Germany}
\email{deicke@uni-paderborn.de}
\author{David Pask}
\author{Iain Raeburn}
\address{Mathematics, University of Newcastle, NSW 2308, Australia}
\email{davidp@maths.newcastle.edu.au, iain@maths.newcastle.edu.au}
\date{December 14, 2001}
\thanks{This research was supported by the Australian Research 
Council and the Deutsche
Forschungsgemeinschaft.}
\maketitle

\section{introduction}

The Cuntz-Krieger algebra $C^*(E)$ of a directed graph $E$ is generated by
a family of partial isometries satisfying relations  which reflect the path
structure of $E$. These \emph{graph algebras} have a rich structure which is
determined by the distribution of loops in the graph. Graph algebras 
have now arisen
in many different situations, and there is increasing interest in 
their interaction
with other graph-theoretic ideas.

Here we consider coverings of directed graphs: morphisms $p:F\to E$ of directed
graphs which are local isomorphisms. We show that the graph algebra 
$C^*(F)$ can be
recovered from $C^*(E)$ as a crossed product by a coaction of a 
homogeneous space
associated to the fundamental group $\pi_1(E)$. The crossed products 
which arise
this way are unusually tractable because we know so much about graph 
algebras, and in
the second half of the paper we give some evidence that these crossed 
products of
graph algebras will be good models for the general theory of crossed products
by homogeneous spaces.

Our results are motivated by work of Kumjian and Pask on \emph{regular
coverings}, which are the orbit maps $p:F\to E$ associated to free 
actions of a group
$G$ on a directed graph $F$ \cite{kp}. They used a description of
$F$ as a skew product
$E\times_c G$ due to Gross and Tucker \cite[Theorem~2.2.2]{gt} to prove that
$C^*(F)\times G$ is stably isomorphic to $C^*(E)$. In \cite{kqr}, it 
was shown that
the description of $F$ as a skew product $E\times_c G$ gives a realisation of
$C^*(F)$ as a crossed product by a coaction $\delta_c$ of $G$ on 
$C^*(E)$, and the
Kumjian-Pask theorem then follows from Katayama duality.

A non-regular covering $p:F\to E$ can also be realised as a kind of skew
product, though the fibre is now the homogeneous space 
$\pi_1(E)/p_*\pi_1(F)$ rather
than a group (see \cite{gt1}, \cite{prho} or \S\ref{sec-cmaps} 
below), and we use this
to obtain our description of $C^*(F)$ as a crossed product 
(Theorem~\ref{provelater}).
It is not in general clear what one should mean by a coaction of a 
homogeneous space
$G/H$ on a
$C^*$-algebra $A$, but here we actually have a normal coaction $\delta$ of the
larger group $G$, and we can use the definition
\begin{equation}\label{dagger}
A\times_\delta (G/H):=\clsp\{j_A(a)j_G(f):a\in A, f\in c_0(G/H)\} \subset
M(A\times_\delta G),
\end{equation}
where $(j_A,j_G)$ is the universal covariant representation of $(A,c_0(G))$ in
$M(A\times_\delta G)$. For non-normal coactions, (\ref{dagger}) is 
the analogue of a
reduced crossed product rather than a full one (see 
Remark~\ref{bloodydefs}), but for
those interested in graph algebras rather than nonabelian duality 
this definition will
suffice.

We begin in Section~\ref{sec-cmaps} by reviewing the structure of 
covering graphs,
and  give a short proof of a recent result of Pask and Rho \cite{prho} which
identifies an arbitrary connected covering as a skew product by a 
homogeneous space.
We prove our main theorems in Section \ref{sec-graph}, and combine them with
results from
\cite{kp} to obtain two corollaries about crossed products of graph algebras,
both of which suggest interesting conjectures about crossed products of
$C^*$-algebras by coactions of a discrete group. In the last two 
sections we prove
these conjectures. In Theorem~\ref{dec_coaction}, we prove that if 
$\delta$ is a
coaction of $G$ on $A$ and $H$ is a subgroup of $G$, then
$A\times_\delta G$ decomposes as an iterated crossed product 
$(A\times (G/H))\times
H$;  the surprise is that we can do this with a crossed product 
rather than a twisted
crossed product. In Theorem~\ref{thm_proper}, we show that the dual 
action of $H$ on
$A\times_\delta G$ is proper in the sense of Rieffel 
\cite{rieffel:90}, and that Rieffel's
generalized fixed-point algebra can be identified with 
$A\times_\delta (G/H)$. This is partial
confirmation of the observation in
\cite{aHRW2} that Rieffel's theory of proper actions, as developed in
\cite{rieff99, M, aHRW}, might provide a useful framework for 
studying crossed products of
$C^*$-algebras by homogeneous spaces.

\section{Coverings of directed grpahs}\label{sec-cmaps}

Let $E=(E^0,E^1,r,s)$ be a directed graph. By a \emph{walk} in $E$ we 
mean a path $a$
in the underlying undirected graph: thus $a=a_1a_2\cdots a_n$, where 
each $a_i$ is
either an edge in $E$ or a reverse edge $e^{-1}$ obtained by 
traversing an edge $e$
backwards. We write $r(e^{-1})=s(e)$, $s(e^{-1})=r(e)$, $r(a)=r(a_n)$ and
$s(a)=s(a_1)$. The \emph{reduction} of a walk $a$ is obtained by deleting any
subwalks of the form
$ee^{-1}$ or $e^{-1}e$; the \emph{reduced product}
$ab$ of two walks with $r(a)=s(b)$ is the reduction of the 
concatenation of $a$ and
$b$.
With this operation, the reduced walks $a$ with $r(a)=s(a)=u$ form a 
group $\pi_1(E,u)$; the
inverse of a reduced walk $a$ is the walk $a^{-1}$ obtained by traversing
$a$ backwards. If
$E$ is connected in the sense that there is a walk between any two 
vertices, then $\pi_1
(E,u)$ is up to isomorphism independent of $u$, and is called the {\em
fundamental group of $E$}. Thus $\pi_1(E,u)$ by definition consists 
of reduced loops based at
$u$. This definition of $\pi_1$ coincides with that in
\cite[2.1.6]{st}, for example, because reducing and path equivalence 
amount to the
same thing.

A surjective morphism $p : F \rightarrow
E$ of directed graphs is a {\em covering} if for each $v \in F^0$,
$p$ maps $r^{-1} (v)$ bijectively onto $ r^{-1} (p ( v ))$ and
$s^{-1} (v)$ bijectively onto $s^{-1} (p ( v ))$ (cf. \cite[2.2.1]{st}). Every
covering $p$ takes reduced walks to reduced walks, respects the 
reduced product and
inverses, and  has the {\em unique walk-lifting property}: if $a$ is 
a reduced walk
with range or source $u\in E^0$, then for every $w \in p^{-1} (u)$ 
there is a unique
reduced walk $\tilde{a}_w$ with range or source $w$ such that $p ( 
\tilde{a}_w ) = a$.
Thus every connected covering $p:F\to E$ induces an injective homomorphism
$p_* : \pi_1 ( F,v) \rightarrow \pi_1 ( E, p(v) )$. The index of 
$p_*\pi_1(F,v)$ is
the number of sheets in the covering:

\begin{lem} \label{theta}
Let $p : F \rightarrow E$ be a connected graph covering
and let $v \in F^0$. For $w \in p^{-1}(p(v))$, choose a reduced walk 
$a$ from $w$ to
$v$. Then $p(a)$ is a reduced loop at $p(v)$, the coset $p(a) p_* 
\pi_1 (F,v)$ is
independent of the choice of
$a$, and the map $\theta : w \mapsto p(a) p_* \pi_1 (F,v)$ is a
bijection of $p^{-1} ( p(v) )$ onto $\pi_1 (E,p(v)) / p_* \pi_1 (F,v)$.
\end{lem}

\begin{proof}
Because $p(w)=p(v)$, $p(a)$ is a loop at $p(v)$, and it is reduced 
because $p$ maps
reduced walks to reduced walks. If
$b$ is any other reduced walk in
$F$ from
$w$ to $v$, then $a^{-1} b$ is a reduced loop based at $v$, so $p 
(a)^{-1} p (b)=p_* (
a^{-1} b )\in p_* \pi_1 (F,v)$ and $p (b) p_* \pi_1 (F,v) = p (a) p_* 
\pi_1 (F,v)$.
Hence $\theta$ is well-defined.

Now suppose $w , w' \in p^{-1} (p(v))$ satisfy $\theta (w) = \theta
(w')$. Let $a$ and $a'$ be reduced walks in $F$ from $w$ and $w'$  to
$v$, so that $\theta (w) = \theta
(w')$ says $p (a) p_* \pi_1 (F,v) = p ( a' ) p_* \pi_1 (F,v)$. Then there
is a reduced walk $d \in \pi_1 (F,v)$ such that $p (a' ) = p (a) p (d) = p
(a d)$. Since $a'$ and $ad$ both terminate at $v$ and satisfy $p(a') 
= p(ad)$, it
follows by unique walk lifting that $w =w'$. Thus $\theta$ is 
injective. If $b \in
\pi_1 (E,p(v))$, let $\tilde{b}_v$ be the lift of $b$ in $F$ with 
range $v$; then $w
:= s (
\tilde{b}_v ) \in p^{-1} (p(v))$ and $\theta ( w ) = p ( \tilde{b}_v 
) p_* \pi_1 (F,v)
= bp_*\pi_1 (F,v)$.
\end{proof}

A \emph{labelling} of a graph $E$ by a group $G$ is a function $c : E^1
\rightarrow G$. If $c$ is a labelling and
$H$ is a subgroup of $G$, the {\em relative skew product} $E \times_{c}
(G/H)$ is the graph  with
$(E\times_{c} (G/H) )^i = E^i \times (G/H)$ for $i=0,1$, $r ( e,gH) = 
( r(e),gH )$ and
$s (e,gH )=(s(e), c (e) gH )$. When $H$ is the trivial subgroup, this 
is the skew
product $E \times_c G$ used in
\cite{kqr} (as opposed to those used in \cite{kp} or \cite[\S 2.3.2]{gt}).

The maps $(x,gH)\mapsto x$ from $(E\times_{c} (G/H) )^i$ to $E^i$ are 
a covering of
$E$, and we are going to prove that every connected covering
$p:F\to E$ has this form. To define a suitable labelling
we choose a
\emph{spanning tree}
$T$ for $E$, which is  a {\em tree} (a connected graph for which 
there is precisely
one reduced walk between any two vertices) with the same vertex set 
as $E$; every
connected directed graph has such a tree (see
\cite[2.1.5]{st}). Now fix a vertex $u\in E^0$. Then for each $w \in 
E^0$, we let
$a_w$ be the unique walk in $T$ from
$u$ to $w$, and define $c=c_{u,T}: E^1 \rightarrow \pi_1
(E,u)$ by
$c(e) = a_{s(e)} e a_{r(e)}^{-1}$.

\begin{prop} \label{cgraphchar}
Let $p : F \rightarrow E$ be a connected covering and $v\in F^0$. Choose a
spanning tree $T$ for $E$, and  define $c:=c_{p(v),T}$. Then $F$ is 
isomorphic to the
skew product $E
\times_{c} (\pi_1 (E,p(v)) / p_* \pi_1 (F,v) )$.
\end{prop}

\begin{proof}
For $z\in F^0$, $a_{p(z)}$ is the unique
reduced walk in $T$ from $u=p(v)$ to $p(z)$. Let $\tilde a_{z}$ be the unique
lifting of $a_{p(z)}$ to a walk in $F$ with $r(\tilde a_{z})=z$, and define
$\tau(z):=s(\tilde a_{z})$, which belongs to $p^{-1}(p(v))$ because
$s(a_{p(z)})=p(v)$. Now define
$\phi : F\to E\times_{c} (\pi_1 (E,p(v)) / p_* \pi_1 (F,v) )$ by
\[
\phi^0 (z) = ( p(z) , \theta ( \tau (z) )  )\  \text{ and }\
\phi^1 ( e ) = ( p(e) , \theta ( \tau ( r(e) ) ) ) .
\]
We shall prove that $\phi$ is the desired isomorphism.

Since $p$ is a graph morphism, it follows immediately from the 
definition of the skew
product that
$r (
\phi (e) )  = \phi ( r(e) )$. Now
\[
\phi (s(e)) = \big( p(s(e)) , \theta ( \tau (s(e)) ) \big)\ \text{ 
whereas }\ s ( \phi (e) ) =
\big( p(s(e)) , c (p(e)) \theta ( \tau (r(e)) ) \big),
\]
so we have to prove that the second coordinates agree.
Recall from Lemma \ref{theta} that $\theta (\tau (s(e))) = p(a) p_* 
\pi_1 (F,v)$ where
$a$ is any reduced walk in $F$ from $\tau (s(e))$ to $v$, and we may choose
$a:=\tilde a_{s(e)} e \tilde a_{r(e)}^{-1} d$ where  $d$ is a reduced walk
from
$\tau ( r(e) )$ to $v$. Then
\begin{align*}
\theta ( \tau (s(e)) )&= p ( \tilde a_{s(e)} e \tilde a_{r(e)}^{-1} d 
) p_* \pi_1
(F,v) \\
&=a_{p(s(e))} p(e) a_{p(r(e))}^{-1} p(d) p_* \pi_1 (F,v) \\
&= c (p(e)) p(d) p_* \pi_1 (F,v)\\
&= c (p(e)) \theta ( \tau (r(e) )).
\end{align*}
Thus $s ( \phi (e) ) = \phi (s(e))$, and $\phi$ is a graph morphism.

To see that
$\phi$ is injective, suppose that $\phi ( e ) = \phi (f)$. Then $p(e) = p(f)$,
and
$\tau ( r(e) ) = \tau ( r(f) )$ because $\theta$ is injective. Since 
$p(e)=p(f)$ we
have
$a_{p(r(e))} = a_{p(r(f))}$, and since $\tau ( r(e) ) = \tau ( r(f) )$ the
liftings $\tilde a_{r(e)}$ and $\tilde a_{r(f)}$ start at the same 
vertex; thus we
can deduce from unique walk-lifting that $\tilde a_{r(e)}=\tilde a_{r(f)}$ and
$r(e)=r(\tilde a_{r(e)})=r(\tilde a_{r(f)})= r(f)$. But now $p(e)=p(f)$ implies
$e=f$ by unique walk-lifting. A similar argument shows that $\phi$ is
injective on $F^0$.

Now suppose that $( e , b p_* \pi_1 (F,v) )$ is an edge
in $E
\times_{c} ( \pi_1 (E,p(v)) / p_* \pi_1 (F,v) )$. Let $\tilde b_v$ be 
the lift of $b$
with range $v$, let $d$ be the lifting of $a_{p(r(e))}$ with 
$s(d)=s(\tilde b_v)$, let
$z=r(d)$, and choose $\tilde e$  to be the unique lift of $e$
with range
$z$. Then by unique walk-lifting, we have $d=\tilde a_z$, so
$\tau(r(\tilde e))$ is by definition $s(\tilde a_{z})=s(\tilde b_v)$, 
and when we
compute $\theta(\tau(r(\tilde e)))$ we can take as $a$ the path 
$\tilde b_v$. Thus
\[
\phi ( \tilde{e} ) = \big( p ( \tilde{e} ) , \theta ( \tau ( r ( 
\tilde{e} ) ) ) \big)
= \big( e, p(\tilde b_v) p_* \pi_1 (F,v) \big) =\big(e,bp_*\pi_1(F,v)\big),
\]
so $\phi$ is surjective on $F^1$. A similar argument shows that 
$\phi$ is surjective
on
$F^0$.
\end{proof}

\begin{rmk}
As it stands, the construction of Proposition~\ref{cgraphchar} 
depends on various
choices, but none of these really matter. Because the covering is connected,
the fundamental group is up to isomorphism independent of the choice of base
point $v$. A different spanning tree $T'$ gives a different labelling $c'$,
but it is cohomologous to $c$ in the sense that there is a function $b:E^0\to
\pi_1(E,p(v))$ such that $b(s(e)) c' (e)  =  c(e) b(r(e))$, and the
corresponding skew products are isomorphic.
\end{rmk}

\section{Graph algebras and coactions}\label{sec-graph}

Our first main theorem realises the $C^*$-algebra of a covering as a 
crossed product
by a coaction. We therefore begin by reviewing our conventions regarding graph
algebras and coactions.

We consider arbitrary directed graphs $E=(E^0,E^1,r,s)$, and 
Cuntz-Krieger $E$-families
consisting of mutually orthogonal projections $\{P_v:v\in E^0\}$ and partial
isometries $\{S_e:e\in E^1\}$ with mutually orthogonal ranges which satisfy
$S_e^*S_e=P_{r(e)}$ and $P_v=\sum_{s(e)=v}S_eS_e^*$ whenever 
$0<|s^{-1}(v)|<\infty$.
The graph algebra $C^*(E)$ is generated by a universal Cuntz-Krieger 
$E$-family $\{p_v,
s_e\}$; for properties of these algebras, see \cite{rsz} and \cite {bhrs}, for
example.

For a discrete group $G$, we denote by $i:G\to UC^*(G)$ the universal unitary
repres\-entation which generates $C^*(G)$. The comultiplication 
$\delta_G:C^*(G)\to
C^*(G)\otimes C^*(G)$ is the integrated form of the representation $s\mapsto
i(s)\otimes i(s)\in U(C^*(G)\otimes C^*(G))$; in this paper, we use 
only the spatial
tensor product of $C^*$-algebras.
A \emph{coaction}
of $G$ on a $C^*$-algebra $A$ is an injective nondegenerate 
homomorphism $\delta:A\to
A\otimes C^*(G)$ such  that
$(\delta\otimes\id)\circ\delta=(\id\otimes\,\delta_G)\circ\delta$. 
We denote by $w_G$
the function $i:G\to C^*(G)$ viewed as an element of
$c_b(G,C^*(G))\subset M(c_0(G)\otimes C^*(G))$. Then a \emph{covariant
representation} of $(A,G,\delta)$ consists of nondegenerate 
representations $\pi:A\to
B(\H)$, $\mu:c_0(G)\to B(\H)$ such that
\[
(\pi\otimes\id)\circ\delta(a)=\mu\otimes\id(w_G)(\pi(a)\otimes
1)\mu\otimes\id(w_G)^*\ \text{ in $M(\K(\H)\otimes C^*(G))$.}
\]
The \emph{crossed product} $A\times_\delta G$ is generated by a 
universal covariant
representation $(j_A,j_G)$ in $M(A\times_\delta G)$, and carries a natural dual
action $\widehat\delta$ of $G$ \cite[\S2]{raeburn:92}.

Because $G$ is discrete, the \emph{spectral subspaces}
\[
A_s:=\{a\in A:\delta(a)=a\otimes i(s)\}
\]
together span a dense subspace of $A$. We write $\chi_t$ for the
characteristic function $\chi_{\{t\}}$ and $a_s$ for a generic 
element of $A_s$; then the
elements
$j_A(a_s)j_G(\chi_t)$ span a dense subspace of $A\times_\delta G$, 
and the dual action is
characterised by
$\widehat{\delta}_r(j_A(a_s)j_G(\chi_t))=j_A(a_s)j_G(\chi_{tr^{-1}})$. 
These facts are
from \cite{quigg:96}, and the following useful lemma is also implicit there.

\begin{lem}\label{quiggcov}
Suppose $\delta$ is a coaction of a discrete group $G$ on a 
$C^*$-algebra $A$, and
$\pi:A\to B(\H)$, $\mu:c_0(G)\to B(\H)$ are nondegenerate representations. Then
$(\pi,\mu)$ is covariant if and only if
\begin{equation}\label{altdefcovrep}
\pi(a_s)\mu(\chi_t)=\mu(\chi_{st})\pi(a_s)\ \text{ for $a_s\in A_s$, $t\in G$.}
\end{equation}
\end{lem}

\begin{proof}
Suppose $(\pi,\mu)$ is covariant, so that
\begin{equation}\label{intermed}
(\pi(a_s)\otimes i(s))\mu\otimes\id(w_G)=\mu\otimes\id(w_G)(\pi(a_s)\otimes 1).
\end{equation}
Slicing $w_G\in c_b(G,C^*(G))$ by the functional $\chi_{st}:z\mapsto z(st)$ on
$C^*(G)$ gives the function $\chi_{st}$ in $c_0(G)$, so applying the slice
map $\id\otimes \chi_{st}$ to (\ref{intermed}) gives
(\ref{altdefcovrep}).

Now suppose (\ref{altdefcovrep}) holds. It suffices to prove 
(\ref{intermed}) for
each element $a_s$ of a spectral subspace, and since $\mu$ is nondegenerate, it
suffices to prove
\begin{equation}\label{3/4}
(\pi(a_s)\otimes
i(s))\mu\otimes\id(w_G)(\mu(\chi_t)\otimes 
1)=\mu\otimes\id(w_G)(\pi(a_s)\otimes
1)(\mu(\chi_t)\otimes 1)
\end{equation}
for every $t\in G$. Using (\ref{altdefcovrep}), the right-hand side 
of (\ref{3/4})
reduces to
\begin{align*}
\mu\otimes\id(w_G)\big((\pi(a_s)\mu(\chi_t))\otimes 1\big)
&=\mu\otimes\id(w_G(\chi_{st}\otimes 1))(\pi(a_s)\otimes 1)\\
&=\mu\otimes\id(\chi_{st}\otimes i(st))(\pi(a_s)\otimes 1)\\
&=\mu(\chi_{st})\pi(a_s)\otimes i(s)i(t)\\
&=(\pi(a_s)\otimes i(s))(\mu(\chi_t)\otimes i(t))\\
&=(\pi(a_s)\otimes i(s))\mu\otimes\id(w_G)(\mu(\chi_t)\otimes 1),
\end{align*}
and the result follows.
\end{proof}

If $H$ is a non-normal subgroup of $G$ there
is some ambiguity about what is meant by the crossed product $A\times_{\delta}
(G/H)$. However, if the coaction $\delta$ is normal in the sense
that there is a covariant representation
$(\pi,\mu)$ of $(A,G,\delta)$ with $\pi$ faithful \cite{q}, then the various
candidates coincide, and we can define
$A\times_{\delta} (G/H)$ using (\ref{dagger}). Since the coactions
in this section are normal, we defer further discussion of this point till
Remark~\ref{bloodydefs}. With this convention, we can state our main theorem.

\begin{thm} \label{provelater}
Let $p : F \rightarrow E$ be a connected covering of directed graphs 
and let $v\in
F^0$. Then there is a normal coaction $\delta$ of $\pi_1 (E,p(v))$ on 
$C^* (E)$ such
that
\[
C^*(F) \cong C^*(E) \times_{\delta} \big(\pi_1(E,p(v))/p_*\pi_1(F,v)\big).
\]
\end{thm}

To prove this theorem, we use  Proposition~\ref{cgraphchar} to 
realise $F$ as the skew
product associated to a labelling $c:E^1\to \pi_1(E,p(v))$. We then recall from
\cite[Lemma~2.3]{kqr} that a labelling $c:E^1\to G$ induces a coaction
$\delta_c:C^*(E)\to C^*(E)\otimes C^*(G)$ characterised by
\begin{equation}\label{chard}
\delta_c(s_e)=s_e\otimes i(c(e))\ \mbox{ and }\ \delta_c(p_v)=p_v\otimes
1_{C^*(G)}.
\end{equation}

\begin{lem} The coaction $\delta_c$ satisfying 
\textnormal{(\ref{chard})} is normal.
\end{lem}

\begin{proof}
Let $\gamma$ denote the gauge action of ${\mathbb T}$ on $C^*(E)$, and choose a
covariant representation $(\pi,U)$ of $(C^*(E),{\mathbb T},\gamma)$ 
such that $\pi$
is faithful. Then $\Ind \pi = ( ( \pi \otimes\lambda )\circ \delta_c 
, 1 \otimes M )$
is a covariant representation of
$(C^*(E),G,\delta_c)$, and
$\{\pi\otimes\lambda(\delta_c(s_e)),\pi \otimes \lambda(\delta_c(p_v))\}$ is a
Cuntz-Krieger $E$-family in which each projection $\pi \otimes
\lambda(\delta_c(p_v))=\pi(p_v)\otimes 1$ is nonzero. Since
\[
(( \pi \otimes \lambda ) \circ \delta_c) ( s_e ) = ( \pi \otimes 
\lambda ) ( s_e
\otimes i(c(e)) ) = \pi ( s_e ) \otimes \lambda_{c(e)},
\]
the representation $U\otimes 1$ implements the gauge action:
\begin{align*}
U_z \otimes 1 ( \pi \otimes \lambda (\delta_c ( s_e )) ) ( U_z \otimes 1 )^*
&= U_z \pi ( s_e ) U_z^* \otimes \lambda_{c(e)}\\
&= \pi ( \gamma_z ( s_e ) ) \otimes \lambda_{c(e)}\\
&=(( \pi \otimes \lambda ) \circ \delta_c) (\gamma_z(s_e)).
\end{align*}
Thus it follows from the gauge-invariant uniqueness theorem 
\cite[Theorem~2.1]{bhrs}
that
$( \pi \otimes \lambda ) \circ \delta_c$ is faithful.\end{proof}

Theorem~\ref{provelater} will therefore follow from 
Proposition~\ref{cgraphchar} and the next
theorem.

\begin{thm} \label{graphsymprod}
Let  $c : E^1 \rightarrow G$ be a labelling of the edges of a 
directed graph $E$ by a
group $G$, and suppose $H$ is a subgroup of $G$. Then
\[
C^* ( E \times_c (G/H) )
\cong C^* (E) \times_{\delta_c} ( G/H).
\]
\end{thm}

\begin{proof} The
covariance relation (\ref{altdefcovrep}) for $(j_A,j_G)$ extends to give
\begin{equation}\label{relationoncosets}
j_A(a_s)j_G(\chi_{tH})=j_G(\chi_{stH})j_A(a_s)\ \mbox{ when $a_s\in A_s$;}
\end{equation}
to see this, just multiply both sides of (\ref{relationoncosets}) by 
$j_G(\chi_r)$, and it
reduces to (\ref{altdefcovrep}).

We now let $\{s_e,p_v\}$ be the canonical Cuntz-Krieger $E$-family generating
$C^*(E)$, and, following \cite[Theorem~2.4]{kqr}, define
\[
t_{f,sH} =j_{C^* (E)} (s_f ) j_{G} ( \chi_{sH} )\ \text{ and
}\       q_{v,tH} = j_{C^*(E)}(p_v)j_G(\chi_{tH}).
\]
We claim that $\{t_{f,sH},q_{v,tH}\}$ is a Cuntz-Krieger $(E\times_c
(G/H))$-family. Because the projections
$p_v$ lie in the spectral subspace
$C^*(E)_{e_G}$ associated to the identity $e_G$ of $G$, 
$j_{C^*(E)}(p_v)$ commutes
with
$j_G(\chi_{tH})$, and hence the $q_{v,tH}$ are mutually orthogonal
projections. Since $s_e\in C^*(E)_{c(e)}$, it follows from 
(\ref{relationoncosets})
that
\begin{align*}
t_{e, sH}^* t_{f,tH} &
= j_G(\chi^*_{sH})j_{C^*(E)}(s_e^* s_f)j_G(\chi_{tH})\\
&=j_{C^*(E)}(s_e^*s_f)j_G(\chi_{c(f)^{-1}c(e)sH}\chi_{tH})\\
&=\begin{cases} j_{C^*(E)}(s_e^* s_f) j_G(\chi_{tH}) &\text{ if 
$c(f)^{-1}c(e)sH=tH$}
\\ 0 & \text{ otherwise} \end{cases}\\
&=\begin{cases} j_{C^*(E)}(p_{r(e)})j_G(\chi_{tH}) &\text{ if $e=f$ 
and $sH = tH$}
\\ 0 & \text{ otherwise} \end{cases}\\
&=\begin{cases} q_{r(e),sH} &\text{ if $(e,sH)=(f,tH)$}
\\ 0 & \text{ otherwise.} \end{cases}
\end{align*}
This proves, first, that the $t_{e,sH}$ are partial isometries with 
initial projections
$q_{r(e,sH)}$, and, second, that the $t_{e,sH}$ have mutually 
orthogonal ranges. A similar
calculation shows that
\[
t_{e,sH}t_{e,sH}^*=j_{C^*(E)}(s_e s_e^*)j_G(\chi_{c(e)sH})\leq
j_{C^*(E)}(p_{s(e)})j_G(\chi_{c(e)sH})=q_{s(e,sH)}.
\]
Since $s^{-1}(v,tH)=\{(e,c(e)^{-1}tH):s(e)=v\}$ is in one-to-one 
correspondence with
$s^{-1}(v)$, we have
$0<|s^{-1}(v,tH)|<\infty$ if and only if
$0<|s^{-1}(v)|<\infty$, and if so,
\begin{align*}
q_{v,tH}&=j_{C^*(E)}(p_v)j_G(\chi_{tH})=\sum_{s(e)=v}j_{C^*(E)}(s_es_e^*)j_G(\chi_{tH})\\
&=\sum_{s(e)=v}j_{C^*(E)}(s_e)j_G(\chi_{c(e)^{-1}tH})j_{C^*(E)}(s_e)^*\\
&=\sum_{s(e,sH)=(v,tH)}t_{e,sH}t_{e,sH}^*.
\end{align*}
Thus $\{t_{f,sH},q_{v,tH}\}$ is a Cuntz-Krieger $(E\times_c
(G/H))$-family, as claimed.  The universal property of the graph 
algebra now gives us
a homomorphism
$\pi_{t,q}$ of
$C^*(E\times_c(G/H))$ into $C^*(E)\times_{\delta_c}(G/H)$ which takes 
the canonical
generating Cuntz-Krieger family to $\{t_{f,sH},q_{v,tH}\}$.

The gauge automorphisms $\gamma_z$
  commute with the coaction
$\delta_c$ in the sense that $\delta_c ( \gamma_z (a) ) = \gamma_z 
\otimes \id(\delta_c (a))$
for $a \in C^* (E)$, and hence by the universal property of  the 
crossed product induce
automorphisms $\gamma_z\times_{\delta_c} G$ of $C^*(E)
\times_{\delta_c} G$. Thus there is
a continuous action $ \gamma \times_{\delta_c} G$ of $\mathbb T$ on $C^*(E)
\times_{\delta_c} G$ such that
\[
( \gamma \times_{\delta_c} G )_z ( j_{C^*(E)}(a)j_G(\chi_{t})) =
j_{C^*(E)}(\gamma_z (a))
j_G(\chi_{t}),
\]
and which leaves the subalgebra $C^* (E) \times_{\delta_c} (G/H) $ of $M(C^*(E)
\times_{\delta_c} G)$ invariant; a calculation shows that 
$\pi_{t,q}\circ\gamma_z
=(\gamma_z\times_{\delta_c} G)\circ \pi_{t,q}$.  Since the coaction
$\delta_c$ is normal, the projections $q_{e,tH}$ are all non-zero, and 
it follows
from the gauge-invariant uniqueness theorem \cite[Theorem~2.1]{bhrs} 
that $\pi_{t,q}$
is injective.

To see that $\pi_{t,q}$ is surjective, it suffices to prove that every
$j_{C^* (E)} ( s_\mu s_\nu^* ) j_{G} ( \chi_{rH} )$ belongs to the range of
$\pi_{t,q}$, or, equivalently, that every $j_{C^* (E)} ( s_\mu ) j_{G} (
\chi_{rH} )j_{C^*(E)}(s_\nu^*)$ with $r(\mu)=r(\nu)$ is in the range 
of $\pi_{t,q}$.
But if
$\tilde\mu$ and
$\tilde\nu$ are the unique liftings of $\mu$ and $\nu$ to paths in 
$E\times_c(G/H)$
with range $(r(\mu),rH)$, then
\[
t_{\tilde\mu}t_{\tilde\nu}^*=j_{C^* (E)} ( s_\mu ) j_{G} (
\chi_{rH} )j_{C^*(E)}(s_\nu^*),
\]
as required.
\end{proof}

\begin{cor} \label{cor1}
Let $\delta_c$ be the coaction of $G$ on $C^*(E)$ induced by a labelling
$c : E^1
\rightarrow G$, and let
$H$ be a subgroup of $G$. Then
$C^* (E) \times_{\delta_c} (G/H)$ is Morita equivalent to $( C^* (E)
\times_{\delta_c} G ) \times_{\widehat{\delta_c}} H$.
\end{cor}

\begin{proof}
The subgroup $H$ acts on the right of $E\times_c G$, and $E \times_c (G/H)$ is
isomorphic to the quotient
$(E\times_c G ) /H$. Thus it follows from Theorem \ref{graphsymprod} that
\[
C^* (E) \times_{\delta_c} (G/H) \cong C^* ( E \times_c (G/H) ) \cong C^* ( ( E
\times_c G ) / H ) .
\]
 From Theorem 1.6 and Corollary 3.1 of \cite{pr} we know
that $C^* ( ( E \times_c G ) / H )$ is Morita equivalent to $C^* ( 
E\times_c G )
\times_\beta H$, where $\beta$ is induced by the right action of $H$ on
$E\times_c G$. To finish off, note that the isomorphism of
$C^*(E\times_c G )$ onto $C^* (E) \times_{\delta_c} G$ carries $\beta$ into the
restriction of the dual action
$\widehat{\delta_c}$.
\end{proof}

\begin{cor} \label{cor2}
Let $\delta_c$ be the coaction of $G$ on $C^*(E)$ induced by a labelling
$c : E^1
\rightarrow G$, and let
$H$ be a subgroup of $G$. Then
there is a coaction $\delta_d$ of $H$ on $C^* (E) \times_{\delta_c} 
(G/H )$ such
that
\[
C^* (E) \times_{\delta_c} G \cong ( C^* (E) \times_{\delta_c}
(G/H)) \times_{\delta_d} H.
\]
\end{cor}

\begin{proof}
  Since $H$ acts freely on $E \times_c G$, the Gross-Tucker theorem
\cite[Theorem~2.2.2]{gt} gives a function $d :(( E \times_c G ) / 
H)^1 \rightarrow H$
such that  $( E \times_c G ) / H \times_d H $ is $H$-isomorphic to 
$E\times_c G$.
Applying Theorem~\ref{graphsymprod} with no subgroup (that is, 
\cite[Theorem~2.4]{kqr}) gives
us a coaction
$\delta_d$ of $H$ such that
\[
C^* ( E \times_c G  ) \cong C^* ( ( E \times_c G ) / H \times_d H ) \cong C^* (
E \times_c (G/H ) ) \times_{\delta_d} H .
\]
We finish off using the ismorphisms
$C^*(E)\times_{\delta_c}G\cong C^*(E\times_c G)$ of \cite[Theorem~2.4]{kqr} and
$C^* ( E \times_c (G/H ) )\cong C^* (E) \times_{\delta_c}
(G/H)$ of Theorem~\ref{graphsymprod}.
\end{proof}

\section{Decomposition of crossed products by discrete coactions}

Our next theorem is motivated by Corollary~\ref{cor2}, which prompts 
the question:
is there always a coaction of $H$ on a crossed product $A\times 
(G/H)$ such that
$A\times G\cong (A\times (G/H))\times H$? We know from 
\cite[Proposition~6.3]{qr}
that this is true when $G$ is discrete and $H$ is normal, because the 
quotient map
$G\to G/H$ certainly admits continuous sections when $G/H$ is 
discrete. Here we show
in rather direct fashion that it is still true for non-normal $H$.

In this section we consider an arbitrary coaction $\delta$ of a 
discrete group $G$ on
a $C^*$-algebra $A$ and a subgroup $H$ of $G$. We now write $A\times_{\delta,r}
(G/H)$ for the $C^*$-algebra
\[
\clsp\{j_A(a)j_G(f):a\in A, f\in c_0(G/H)\} \subset M(A\times_\delta G),
\]
and call it the
\emph{reduced crossed product of $A$ by the homogeneous space $G/H$}.

\begin{rmk}\label{bloodydefs}
Our reasons for calling this a reduced crossed product are a little 
delicate. When
$H$ is a normal subgroup, it makes sense to restrict a coaction $\delta:A\to
A\otimes C^*(G)$ to a coaction $\delta|$ of
$G/H$ and form the usual crossed product. If $\pi$ is a faithful 
representation of
$A$ on $\H$, then the spatially-defined crossed product generated by the operators
\[
(\pi\otimes\lambda^{G/H})\circ \delta|(a), 1\otimes M(f)\in 
B(\H\otimes \l^2(G/H))
\]
has the universal property which characterises $A\times_{\delta|} (G/H)$; since
$(\pi\otimes\lambda^{G/H})\circ \delta|, 1\otimes M)$ is the regular 
representation of
$(A,G/H,\delta|)$ induced from $\pi$, this can be interpreted as 
saying that the full
and reduced crossed products coincide. However, the system 
$(A,G/H,\delta|)$ also has
a natural representation $((\pi\otimes \lambda^G)\circ \delta, 1\otimes M)$ on
$\H\otimes \l^2(G)$, where now $M$ is the representation of $c_0(G)$ 
restricted to
$c_0(G/H)\subset M(c_0(G))$, and this representation on
$\H\otimes \l^2(G)$ makes sense even when $H$ is not normal. So in
\cite{ekr}, the reduced crossed product $A\times_{\delta,r} (G/H)$ 
was by definition
the $C^*$-algebra on $\H\otimes \l^2(G)$ generated by the operators
$(\pi\otimes \lambda^G)\circ \delta(a)(1\otimes M(f))$ for $f\in 
c_0(G/H)$. Since
the representation $((\pi\otimes \lambda^G)\circ \delta)\times (1\otimes M)$ is
faithful on $A\times_\delta G$ and everything lies in 
$M(A\times_\delta G)$, we can
view this $A\times_{\delta, r}(G/H)$ as a subalgebra of 
$M(A\times_\delta G)$. Thus
our new definition of $A\times_{\delta, r}(G/H)$ coincides with that 
of \cite{ekr}.

The notation in \cite{ekr} was consistent with
that of Mansfield \cite{man}, who for normal $H$ used
$A\times_\delta (G/H)$ to distinguish the subalgebra of $B(\H\otimes 
\l^2(G))$  from
the usual crossed product $A\times_{\delta|}(G/H)$  on $\H\otimes \l^2(G/H)$.
Mansfield proved that when $H$ is normal and amenable
$A\times_{\delta|}(G/H)$ is isomorphic to $A\times_{\delta,r}(G/H)$
\cite[Proposition~7]{man}. For nonamenable normal subgroups, though, 
the two algebras
need not be isomorphic --- for example, if $H=G$ and the coaction is 
not normal in
the sense of Quigg \cite{q}.

The Fell-bundle approach gives another way of defining a reduced 
crossed product. If
$\delta:A\to A\otimes C^*(G)$ is a coaction of a discrete group, then 
the  spectral
subspaces
$\{A_s:s\in G\}$ form a Fell bundle $\mathcal{A}$ over $G$. When $H$ 
is a subgroup of
$G$, we can form a Fell bundle $\mathcal{A}\times (G/H)$ over the 
transformation
groupoid $G\times (G/H)$ (see \cite[\S2]{eq}), and this has a regular 
representation as
adjointable operators on a Hilbert module $L^2(\mathcal{A}\times 
G/H)$ which gives a
reduced cross-sectional algebra $C_r^*(\mathcal{A}\times (G/H))$. It 
is proved in
\cite[Proposition~2.10]{eq} that  $C_r^*(\mathcal{A}\times (G/H))$ is naturally
isomorphic to the $C^*$-subalgebra of $M(C^*(\mathcal{A})\times_{\delta^m}G)$
generated by the images of $C^*(\mathcal{A})$ and $c_0(G/H)$, where 
$\delta^m$ is the
maximal coaction on the full cross-sectional algebra 
$C^*(\mathcal{A})$. However, we
know from \cite[Lemma~2.1]{eq0} that there is an isomorphism of
$C^*(\mathcal{A})\times_{\delta^m}G$ onto $A\times_\delta G$ taking
$j_{C^*(\mathcal{A})}(C^*(\mathcal{A}))$ to $j_A(A)$ and matching up 
the copies of
$c_0(G/H)$, so our $A\times_{\delta,r}(G/H)$ is
isomorphic to $C_r^*(\mathcal{A}\times (G/H))$.
\end{rmk}

\begin{thm}\label{dec_coaction}
Let $\delta$ be a coaction of a discrete group $G$ on a $C^*$-algebra 
$A$, let $H$ be
a subgroup of $G$, and let $\sigma:G/H\to G$ be a section such that 
$\sigma(H)=e$.
Then there is a normal coaction $\delta^\sigma$ of $H$ on 
$A\times_{\delta,r}(G/H)$
such that
\begin{equation}\label{eq:delta_s}
\delta^\sigma(j_A(a_r)j_G(\chi_{tH}))
=j_A(a_r)j_G(\chi_{tH})\otimes i(\sigma(rtH)^{-1}r\sigma(tH))
\end{equation}
for $a_r\in A_r$ and $tH\in G/H$, and there is an isomorphism of
$(A\times_{\delta,r}(G/H))\times_{\delta^\sigma}H$ onto 
$A\times_\delta G$ which
carries the dual action $\widehat{\delta^\sigma}$ into the 
restriction of the dual
action~$\widehat \delta$.
\end{thm}

\begin{proof}
Define $\phi:G\to H$ by $\phi(s)=\sigma(sH)^{-1}s$, and let $\phi^*$ denote the
homomorphism $f\mapsto f\circ \phi$ of $c_0(H)$ into 
$c_b(G)=M(c_0(G))$; note that
$\phi^*$ is nondegenerate. The function $w_H$ which takes $h$ to its
canonical image $i(h)$ in $UC^*(H)$ is a unitary element of $c_b(H, 
C^*(H))\subset
M(c_0(H)\otimes C^*(H))$, and is a corepresentation of
$H$.  Thus the map $\delta_1:A\times_\delta G\to(A\times_\delta
G)\otimes C^*(H)$ defined by
\begin{equation}\label{covofdeltac}
\delta_1(c):=((j_G\circ\phi^*)\otimes\operatorname{id})(w_H)(c\otimes 1)
((j_G\circ\phi^*)\otimes\operatorname{id})(w_H^*)
\end{equation}
is a coaction of $H$ on $A\times_\delta G$. Since $\delta_1$ is
nondegenerate as a homomorphism, it extends to $M(A\times_\delta G)$. 
We claim that
this extension of $\delta_1$ maps $A\times_{\delta,r} (G/H)$ into 
$(A\times_{\delta,r}
(G/H))\otimes C^*(H)$, and therefore restricts to a coaction 
$\delta^\sigma$ of $H$
on
$A\times_{\delta,r}(G/H)$.

For $s\in G$, we have $(\chi_s\otimes 1)(\phi^*\otimes \id)(w_H)=\chi_s\otimes
i(\phi(s))$. Thus for $a_r\in A_r$ and $tH\in G/H$ we have
\begin{align*}
(j_G(\chi_s)&\otimes 1)\delta_1(j_A(a_r)j_G(\chi_{tH}))\\
&=\big(j_G(\chi_s)j_A(a_r)j_G(\chi_{tH})\otimes i(\phi(s))\big)
((j_G\circ\phi^*)\otimes\id)(w_H^*)\\
&=\begin{cases}
\big(j_A(a_r)j_G(\chi_{r^{-1}s})\otimes i(\phi(s))\big)
((j_G\circ\phi^*)\otimes\id)(w_H^*)&\text{if $r^{-1}s\in tH$}\\
0&\text{if $r^{-1}s\notin tH$}
\end{cases}\\
&=\begin{cases}
j_A(a_r)j_G(\chi_{r^{-1}s})\otimes i(\phi(s)\phi(r^{-1}s)^{-1})
&\text{if $r^{-1}s\in tH$}\\
0&\text{if $r^{-1}s\notin tH$}
\end{cases}\\
&=(j_G(\chi_s)\otimes 1)\big(j_A(a_r)j_G(\chi_{tH})\otimes
i(\sigma(rtH)^{-1}r\sigma(tH))\big).
\end{align*}
This implies that $\delta_1$ maps
$A\times_{\delta,r} (G/H)$ into $(A\times_{\delta,r}(G/H))\otimes 
C^*(H)$, and that
the restriction $\delta^\sigma:=\delta_1|_{A\times(G/H)}$ satisfies
(\ref{eq:delta_s}). Since $A_e$ contains an approximate identity for
$A$ \cite[Corollary~1.6]{quigg:96}, (\ref{eq:delta_s}) implies that 
$\delta^\sigma$ is
a nondegenerate homomorphism. It follows from the corresponding properties of
$\delta_1$ that $\delta^\sigma$ is injective and satisfies the 
coaction identity. If
$\pi\times \mu$ is a faithful representation of $A\times_\delta G$, 
then $\pi\times
\mu$ is faithful on $A\times_{\delta,r} (G/H)$. In the regular representation
$\big(((\pi\times\mu)\otimes\lambda^H)\circ\delta^\sigma, 1\otimes M\big)$ on
$\H\otimes\l^2(H)$, the algebra component
$((\pi\times\mu)\otimes\lambda^H)\circ\delta^\sigma$ is implemented 
by conjugating
the representation $\big((\pi\times\mu)|_{A\times (G/H)}\big)\otimes 1$ by the
unitary $((\mu\circ\phi^*)\otimes\lambda^H)(w_H)$, and hence is injective. Thus
$\delta^\sigma$ is normal.

To obtain the isomorphism in the Theorem, we prove that if $\iota$ is 
the inclusion
of $A\times_{\delta,r}(G/H)$ in $M(A\times_\delta G)$, then 
$(A\times_\delta G,\iota,
j_G\circ\phi^*)$ is a crossed product for
$(A\times_{\delta,r}(G/H),H,\delta^\sigma)$, in the sense that 
$(A\times_\delta G,\iota,
j_G\circ\phi^*)$ has properties (a), (b) and (c) of
\cite[Definition~2.8]{raeburn:92}.

Equation~(\ref{covofdeltac}) implies that $(\iota,j_G\circ \phi^*)$ 
is covariant for
$\delta^\sigma:=\delta_1|_{A\times(G/H)}$, so (a) holds.  To verify (c),
note first that
\begin{align*}
\big(\chi_{tH}\phi^*(\chi_{\phi(t)})\big)(s)=1
&\Longleftrightarrow \text{$s\in tH$ and $\phi(s)=\phi(t)$}\\
&\Longleftrightarrow \text{$sH= tH$ and $\sigma(sH)^{-1}s=\sigma(tH)^{-1}t$}\\
&\Longleftrightarrow s=t,
\end{align*}
so that $\chi_{tH}\phi^*(\chi_{\phi(t)})=\chi_{t}$. Thus each spanning element
$j_A(a_r)j_G(\chi_t)$ of $A\times_\delta G$ can be written as
\[
j_A(a_r)j_G(\chi_t)=j_A(a_r)j_G(\chi_{tH})j_G(\phi^*(\chi_{\phi(t)}))
=\iota\big(j_A(a_r)j_G(\chi_{tH})\big)j_G\circ \phi^*(\chi_{\phi(t)}),
\]
which implies that the elements $\iota(b)j_G\circ \phi^*(\chi_{h})$ span a
dense subspace of $A\times_\delta G$.

It remains to establish \cite[Definition~2.8(b)]{raeburn:92}. Suppose
$(\rho,\mu)$ is a covariant representation of
$(A\times_{\delta,r}(G/H),H,\delta^\sigma)$. Both
$j_A(A)$ and $j_G(c_0(G/H))$ multiply $A\times_{\delta,r}(G/H)$, and 
hence we can
view $j_A$ and $j_G$ as nondegenerate homomorphisms into
$M(A\times_{\delta,r}(G/H))$. By
\cite[Proposition~2.6]{eq}, we have
\[
j_A(a_r)j_G(\chi_{tH})=j_G(\chi_{stH})j_A(a_r)\ \mbox{ for $a_r\in A_r$}.
\]
Now let $\pi:=\rho\circ j_A$ and $\nu:=\rho\circ j_G|_{c_0(G/H)}$. 
Then for $a_e\in
A_e$, $j_A(a_e)j_G(\chi_{tH})$ belongs to the spectral
subspace $(A\times_{\delta,r}(G/H))_e$, and the covariance of 
$(\rho,\mu)$ implies
that
$\pi(a_e)\nu(\chi_{tH})=\rho(j_A(a_e)j_G(\chi_{tH}))$ commutes with 
$\mu(\chi_h)$ for
$h\in H$. By letting $a_e$ run through an approximate identity in 
$A_e$, we conclude
that $\nu(\chi_{tH})$ and $\mu(\chi_h)$ commute. Thus they combine to give a
representation $\nu\otimes\mu$ of $c_0(G/H)\otimes c_0(H)$. The 
section $\sigma$
gives an isomorphism $\psi$ of $c_0(G/H)\otimes c_0(H)$ onto $c_0(G)$ such that
$\psi(f\otimes g)(s)=f(sH)g(\phi(s))$;  pulling $\nu\otimes\mu$ over gives a
representation $\omega:=(\nu\otimes\mu)\circ\psi^{-1}$ of $c_0(G)$ satisfying
$\omega(\chi_t)=\nu(\chi_{tH})\mu(\chi_{\sigma(tH)^{-1}t})$ for $t\in 
G$. We claim
that $(\pi,\omega)$ is a covariant representation of
$(A,G,\delta)$. To see this, let $a_r\in A_r$  and $t\in G$. Then
\begin{align*}
\pi(a_r)\omega(\chi_t)&=\rho(j_A(a_r)j_G(\chi_{tH}))\mu(\chi_{\sigma(tH)^{-1}t})\\
&=\mu(\chi_{\sigma(rtH)^{-1}r\sigma(tH)\sigma(tH)^{-1}t})\rho(j_A(a_r)j_G(\chi_{tH}))
\quad\mbox{ by (\ref{eq:delta_s})}\\
&=\mu(\chi_{\sigma(rtH)^{-1}rt})\pi(a_r)\nu(\chi_{tH})\\
&=\mu(\chi_{\sigma(rtH)^{-1}rt})\nu(\chi_{rtH})\pi(a_r)\\
&=\omega(\chi_{rt})\pi(a_r).
\end{align*}
It follows from Lemma~\ref{quiggcov} that
$(\pi,\omega)$ is covariant, and hence gives a nondegenerate 
representation $\pi\times
\omega$ of $A\times_\delta G$. This extends to $M(A\times_\delta G)$, 
and satisfies
\[
(\pi\times
\omega)(j_A(a_r)j_G(\chi_{tH}))=\pi(a_r)\omega(\chi_{tH})=\pi(a_r)\nu(\chi_{tH})=
\rho(j_A(a_r)j_G(\chi_{tH}));
\]
since $\psi(1\otimes \chi_h)=\phi^*(\chi_h)$, we also have
\[
(\pi\times
\omega)(j_G\circ\phi^*(\chi_h))=
\omega(\phi^*(\chi_h))=\nu\otimes\mu(1\otimes\chi_h)=\mu(\chi_h).
\]

This completes the proof that
$(A\times_\delta G,\iota,j_G\circ \phi^*)$ is a crossed product for the
coaction $\delta^\sigma$ of $H$ on
$A\times_{\delta,r} (G/H)$. There is therefore an isomorphism
$\theta$ of $A\times_\delta G$ onto
$(A\times_{\delta,r} (G/H))\times_{\delta^\sigma} H$ such that
$\theta\circ\iota=j_{A\times(G/H)}$ and $\theta\circ j_G\circ \phi^*=j_H$.
It remains to prove that $\theta$ converts $\widehat{\delta}|_H$ into
$\widehat{\delta^\sigma}$. But for  $k\in H$, both $\widehat{\delta}_k$ and
$(\widehat{\delta^\sigma})_k$ fix the copies of $A$ and are given by right
translation $\rt_k$ on $c_0(H)$ and $c_0(G)$, respectively.  Since
$\chi_{tH}\phi^*(\chi_h)=\chi_{\sigma(tH)h}$, we have
\begin{align*}
(\theta\circ\widehat{\delta}_k)\big(j_A(a_r)j_G(\chi_{tH})j_G({\phi}^*(\chi_h))\big)
&=\theta\big(\widehat{\delta}_k(j_A(a_r)j_G(\chi_{\sigma(tH)h}))\big)\\
&=\theta\big(j_A(a_r)j_G(\chi_{\sigma(tH)hk^{-1}})\big)\\
&=\theta\big(j_A(a_r)j_G(\chi_{tH})j_G({\phi}^*(\chi_{hk^{-1}}))\big)\\
&=j_A(a_r)j_G(\chi_{tH})j_H(\chi_{hk^{-1}})\\
&=(\widehat{\delta^\sigma})_k\big(j_A(a_r)j_G(\chi_{tH})j_H(\chi_h)\big)\\
&=((\widehat{\delta^\sigma})_k\circ\theta)\big(j_A(a_r)j_G(\chi_{tH})j_G({\phi}^*(\chi_h))\big),
\end{align*}
as required.
\end{proof}

\begin{cor} \label{cor_dec_coaction}
Let $(A,G,\delta)$ and $H$ be as in Theorem \ref{thm_proper}. Then
\begin{equation}\label{nonstableiso}
(A\times_{\delta,r}
(G/H))\otimes\mathcal{K}(L^2(H))\cong (A\times_\delta 
G)\times_{\widehat{\delta}|_H,r}
H.
\end{equation}
\end{cor}

\begin{proof}
It follows from Theorem~\ref{dec_coaction} that
\[
((A\times_{\delta,r} 
(G/H))\times_{\delta^\sigma}H)\times_{\widehat{\delta^\sigma},r}H
\cong (A\times_\delta G)\times_{\widehat{\delta}|_H,r}H.
\]
Since $\delta^\sigma$ is normal, it follows from 
\cite[Theorem~2.10]{quigg:91} that
Katayama duality  works:
\[
((A\times_{\delta,r} 
(G/H))\times_{\delta^\sigma}H)\times_{\widehat{\delta^\sigma},r}H
\cong (A\times_{\delta,r} (G/H))\otimes\mathcal{K}(L^2(H)).
\]
Composing the isomorphisms gives the result.
\end{proof}

\begin{rmks}
(a) Corollary~\ref{cor_dec_coaction} is slightly stronger than
\cite[Theorem 3.4]{eq}, which says that
$A\times_{\delta,r}(G/H)$ is Morita equivalent to
$(A\times_\delta G)\times_{\widehat{\delta}|_H,r}H$. The argument in 
\cite{eq} uses
Fell-bundle techniques to build an imprimitivity bimodule which 
implements the Morita
equivalence. We can get a third proof of  \cite[Theorem 3.4]{eq} by 
combining our
Theorem~\ref{thm_proper} below and \cite[Corollary~1.7]{rieffel:90}.

(b) There is an analogous Morita equivalence for full crossed products
in \cite[Theorem 3.1]{eq}, and it is tempting to ask whether our 
methods will also
give an isomorphism like (\ref{nonstableiso}) for full crossed 
products. For the full
crossed product
$A\times_{\delta,f}(G/H)$ we take the full cross-sectional algebra
$C^*(\mathcal{A}\times (G/H))$ of the Fell bundle $\mathcal{A}\times 
(G/H)$ whose
reduced cross-sectional algebra $C^*_r(\mathcal{A}\times
(G/H))$ we have already seen to be isomorphic to 
$A\times_{\delta,r}(G/H)$. We can use the
universal property of $C^*(\mathcal{A}\times (G/H))$ to see that
there is a coaction $\delta^{\sigma,f}$ of $H$ on
$A\times_{\delta,f}(G/H)$ satisfying (\ref{eq:delta_s}), and then the 
arguments of
Theorem~\ref{dec_coaction} show that we have
\[
\big((A\times_{\delta,f} (G/H))
\times_{\delta^{\sigma,f}}G,H,\widehat{\delta^{\sigma,f}}\big)\cong
(A\times_\delta G,H,\widehat{\delta}|_H).
\]
If we knew that $\delta^{\sigma,f}$ were \emph{maximal} in the sense 
of \cite{ekr},
so that
\[
(A\times_{\delta,f}(G/H))\otimes\mathcal{K}(L^2(H))\cong 
\big((A\times_{\delta,f}
(G/H))\times_{\delta^\sigma}H\big)\times_{\widehat{\delta^\sigma}}H,
\]
then we could deduce a version of Corollary~\ref{cor_dec_coaction} 
for full crossed
products. We have, however, been unable to prove that 
$\delta^{\sigma,f}$ is maximal.
\end{rmks}

\section{The crossed product as a generalized fixed-point algebra}

Our last result is suggested by Corollary~\ref{cor1} and
Corollary~\ref{cor_dec_coaction}. We expect crossed products by free actions to
be related to the fixed-point algebra for the action. Thus these 
corollaries suggest
that we might be able to view $A\times_{\delta,r}(G/H)$ as a 
fixed-point algebra for
the dual action $\widehat{\delta_c}$. We shall confirm this by proving that
$\widehat{\delta_c}$ is proper and saturated in the sense of Rieffel
\cite{rieffel:90}, so that there is a generalized fixed-point algebra 
which is Morita
equivalent to the reduced crossed product, and that this fixed-point algebra is
precisely
$A\times_{\delta,r}(G/H)$.

We begin by recalling Rieffel's definition of proper action as it applies to
an action $\alpha$ of a discrete group $G$ \cite{rieffel:90}. He says 
that $\alpha$
is \emph{proper} if there is a dense $\alpha$-invariant
$*$-subalgebra $A_0$ of $A$ such that
\begin{enumerate}
\item[(i)] for every $a,b\in A_0$, the function ${}_E\langle a,b 
\rangle : s \mapsto
a \alpha_s (b^*)$ is in $l^1(G,A)$, and
\item[(ii)] for every $a, b\in A_0$, there is
a multiplier $\langle a,b \rangle_D$ of $A$ such that  $A_0\langle a,b
\rangle_D\subset A_0$, $\langle a,b \rangle_D A_0\subset A_0$, and
\begin{equation}\label{eq:proper}
\sum_{s \in G} c \alpha_s ( a^*b ) =c \langle a,b \rangle_D \ \mbox{ 
for all $c\in
A_0$.}
\end{equation}
(The sum on the left converges absolutely by (i), so 
(\ref{eq:proper}) makes sense;
the multipliers $\langle a,b
\rangle_D$ satisfying (\ref{eq:proper}) are automatically invariant 
under $\alpha$.)
\end{enumerate}
If $\alpha$ is proper the {\em generalized fixed-point algebra}
$A^\alpha:=\clsp\{\langle a,b\rangle_D:a,b\in A_0\}$ is a $C^*$-subalgebra of
$M(A)$. The action
$\alpha$ is {\em saturated} if $\newspan\{{}_E\langle a,b
\rangle:a,b \in A_0\}$ is dense in the reduced crossed product 
$A\times_{\alpha,r}
G$.

Corollary~1.7 of \cite{rieffel:90} says that if $\alpha$ is proper 
and saturated,
then $A\times_{\alpha,r}G$ is Morita equivalent to $A^\alpha$.

\begin{thm}\label{thm_proper}
   Let $\delta$ be a coaction of a group $G$ on a $C^*$-algebra $A$, 
and let $H$ be a
subgroup of
   $G$. Then the dual action  $\widehat{\delta}|_H$ of $H$ on $A\times_\delta G$
   is proper and saturated, with generalized fixed-point algebra
\begin{equation}\label{idfpa}
(A\times_\delta
   G)^{\widehat{\delta}|_H}=
  \clsp \{j_A(a_r)j_G(\chi_{tH}):r,t\in G,a_r\in A_r\}=A\times_{\delta,r}(G/H).
\end{equation}
\end{thm}

\begin{proof}
For the dense $\widehat{\delta}|_H$-invariant $\ast$-subalgebra we take
\[
(A\times_\delta G)_0:=\operatorname{span} 
\{j_A(a_r)j_G(\chi_{t}):r,t\in G, a_r\in
A_r\},
\]
which is dense because $A=\clsp\{a_r\}$ and $c_0(G)=\clsp\{\chi_t\}$, a
$*$-subalgebra because
$j_A(a_r)j_G(\chi_t)=j_G(\chi_{rt})j_A(a_r)$, and invariant because
$\widehat{\delta}_h(j_A(a_r)j_G(\chi_t))=j_A(a_r)j_G(\chi_{th^{-1}})$.

To verify (i), let $a=j_A(a_r)j_G(\chi_t)$, $b=j_A(b_s)j_G(\chi_u)$, 
and $h\in H$.
Then
\begin{align}\label{eq:inner_product_E}
{}_E\langle
a,b\rangle(h) 
&=j_A(a_r)j_G(\chi_t)\widehat{\delta}_h((j_A(b_s)j_G(\chi_u))^*) \\
&=j_A(a_r)j_G(\chi_t)\widehat{\delta}_h(j_G(\chi_u)
j_A(b_s)^*)\notag\\
&=j_A(a_r)j_G(\chi_t)j_G(\chi_{uh^{-1}})j_A(b_s)^*\notag\\
&=\begin{cases}
j_A(a_r)j_G(\chi_t)j_A(b_s)^*&\text{if $h=t^{-1}u$}\\
0&\text{otherwise,}
\end{cases}\notag
\end{align}
so the support of $\langle a,b\rangle_E$ consists of the single point 
$t^{-1}u$. In
general $\supp\langle a,b\rangle_E$ is finite for $a,b\in 
(A\times_\delta G)_0$, so
${}_E\langle a,b\rangle$ is certainly summable. Thus (i) holds.

To show (ii), we define
\begin{equation}\label{Dip}
\langle j_A( a_r ) j_G( \chi_t ), j_A( b_s) j_G( \chi_u ) \rangle_D:=
\begin{cases}
      j_A( a_r^* b_s ) j_G ( \chi_{uH}) & \text{if $rt=su$} \\
      0 & \text{otherwise;}
     \end{cases}
\end{equation}
once we have shown that these have the required property 
(\ref{eq:proper}), we can
find suitable multipliers $\langle a,b\rangle_D$ for general
$a,b\in (A\times_\delta G)_0$ by adding up multipliers of the form (\ref{Dip}).
Meanwhile, write $a=j_A( a_r ) j_G( \chi_t )$, $b=j_A( b_s) j_G( \chi_u )$ and
$c=j_A(c_{\ell})j_G(\chi_v)$. Then the formula
\[
c\langle a,b\rangle_D=
\begin{cases}
      j_A(c_{\ell}a_r^*b_s)j_G(\chi_{s^{-1}rv}\chi_{uH}) & \text{if $rt=su$} \\
      0 & \text{otherwise}
     \end{cases}
\]
and a similar one for $\langle a,b\rangle_Dc$ show that $c\langle 
a,b\rangle_D$  and
$\langle a,b\rangle_Dc$ belong to
$(A\times_\delta G)_0$.

We now verify (\ref{eq:proper}) for $a$, $b$ and $c$ as above. Both sides of
(\ref{eq:proper}) vanish unless $rt=su$. If $rt=su$ and $h\in H$, then
\begin{align*}
c\widehat{\delta}_h(a^*b)&=j_A(c_\ell
a_r^*b_s)j_G(\chi_{s^{-1}rv}\chi_{s^{-1}rth^{-1}})\\
&=\begin{cases}
j_A(c_\ell a_r^*b_s)j_G(\chi_{uh^{-1}})&\text{if $h=v^{-1}t$}\\
0&\text{otherwise,}
\end{cases}
\end{align*}
so
$\sum_{h\in H}c\widehat{\delta}_h(a^*b)$ is nonzero only if 
$v^{-1}t\in H$, and then
\begin{align*}
\sum_{h\in H}c\widehat{\delta}_h(a^*b)
&=j_A(c_\ell a_r^*b_s)j_G(\chi_{ut^{-1}v})\\
&=j_A(c_\ell a_r^*b_s)j_G(\chi_{s^{-1}rv})\\
&=j_A(c_\ell a_r^*b_s)j_G(\chi_{s^{-1}rv}\chi_{uH})\ \mbox{ because
$u^{-1}s^{-1}rv=(v^{-1}t)^{-1}\in H$}\\
&= c\langle a,b\rangle_D.
\end{align*}
Thus $\widehat{\delta}|_H$ is proper.

Since the spectral subspace $A_e$ contains an approximate identity for $A$
\cite[Corollary 1.6]{quigg:96}, taking $s=e$ in
(\ref{eq:inner_product_E}) shows that the elements of the form
${}_E\langle a, b\rangle$ are dense in $\ell^1(H, A\times_\delta G)$, hence in
$(A\times_\delta G)\times_{\widehat\delta, r} H$. Thus
$\widehat{\delta}|_H$ is saturated.

It remains to establish (\ref{idfpa}). The formula (\ref{Dip}) shows that
$(A\times_\delta G)^{\widehat{\delta}|_H}$ is contained in the
right-hand side. Let  $a_r\in A_r$ and $tH\in G/H$. Then if 
$(u_\lambda)_{\lambda}$ is
an approximate identity for $A$ contained in $A_e$, we have
\[
j_A(a_r)j_G(\chi_{tH})=\lim_\lambda j_A(u_\lambda a_r)j_G(\chi_{tH})
=\lim_\lambda\langle j_A(u_\lambda) j_G(\chi_{st}),
j_A(a_r)j_G(\chi_t)\rangle_D,
\]
which is in $(A\times_\delta G)^{\widehat{\delta}|_H}$.
\end{proof}


\begin{thebibliography}{20}

\bibitem{bhrs} T.~Bates, J.H.~Hong, I.~Raeburn and W. Szyma{\' n}ski.
\newblock The ideal structure of $C^*$-algebras of infinite graphs,
\newblock preprint, Univ. of Newcastle, 2001.

\bibitem{ekr} S. Echterhoff, S. Kaliszewski and J. Quigg. Maximal 
coactions, preprint,
Arizona State Univ., 2001.

\bibitem{ekr} S. Echterhoff, S. Kaliszewski and I. Raeburn. Crossed 
products by dual
coactions of groups and homogeneous spaces, {\em J. Operator Theory 
{\bf 39}} (1998),
151--176.

\bibitem{eq0} S~Echterhoff and J. Quigg. Induced coactions of 
discrete groups on
$C^*$-algebras, {\em Canad. J. Math. {\bf 51}} (1999), 745--770.

\bibitem{eq} S~Echterhoff and J. Quigg.
\newblock Full duality for coactions of discrete groups,
\newblock {\em Math. Scand.}, to appear.

\bibitem{aHRW} A. an Huef, I. Raeburn and D.P. Williams. Proper 
actions on imprimitivity
bimodules and decompositions of Morita equivalences, preprint, Univ. 
of Newcastle, 2001.

\bibitem{aHRW2} A. an Huef, I. Raeburn and D.P. Williams. A symmetric 
imprimitivity theorem
for commuting proper actions, in preparation.



\bibitem{gt1} J.L.~Gross and T.W.~Tucker.
\newblock Generating all graph coverings by permutation voltage assignments,
\newblock {\em Discrete Math. {\bf 18}} (1977), 273--283.

\bibitem{gt} J.L.~Gross and T.W.~Tucker.
\newblock {\em Topological Graph Theory,}
\newblock Wiley Interscience Series in Discrete Math. \& Optimization,
1987.

\bibitem{kqr} S.~Kaliszewski, J.~Quigg and I.~Raeburn.
\newblock Skew products and crossed products by coactions,
\newblock {\em J. Operator Theory}, to appear.

\bibitem{kp} A.~Kumjian and D.~Pask.
\newblock $C^*$-algebras of directed graphs and group actions,
\newblock {\em Ergod. Th. \& Dynam. Sys. {\bf 19}} (1999), 1503--1519.

\bibitem{man} K. Mansfield. Induced representations of crossed 
products by coactions,
{\em J. Funct. Anal. {\bf 97}} (1991), 112--161.

\bibitem{M} R. Meyer. Generalized fixed-point algebras and 
square-integrable group actions,
{\em J. Funct. Anal.}, to appear.

\bibitem{pr} D.~Pask and I.~Raeburn.
\newblock Symmetric imprimitivity theorems for graph $C^*$-algebras,
\newblock {\em Internat. J. Math. {\bf 12}} (2001), 609--623.

\bibitem{prho} D.~Pask and S.-J.~Rho.
\newblock Intrinsic properties of simple graph algebras,
\newblock {\em Proceedings of OAMP Conference,  Constanza, 2001}, to appear.

\bibitem{quigg:91}
J. Quigg. Full $C^*$-crossed product duality, {\em J. Austral. Math. 
Soc. (Series A)
{\bf {50}}} (1991), 34--52.

\bibitem{q} J. Quigg.
\newblock Full and reduced $C^*$-coactions,
\newblock {\em Math. Proc. Camb. Phil. Soc. {\bf 116}} (1995), 435--450.

\bibitem{quigg:96} J.~Quigg.
\newblock Discrete $C^*$-coactions and $C^*$-algebraic bundles,
\newblock {\em J. Austral. Math. Soc. (Series A) {\bf {60}}} (1996), 204--221.



\bibitem{qr} J. Quigg and I. Raeburn.
\newblock Induced $C^*$-algebras and Landstad
duality for twisted coactions,
\newblock {\em Trans. Amer. Math. Soc. {\bf 347}}
(1995), 2885--2915.

\bibitem{raeburn:92} I.~Raeburn.
\newblock On crossed products by coactions and their representation theory,
\newblock {\em Proc. London Math. Soc. \bf {64}} (1992), 625--652.

\bibitem{rsz} I. Raeburn and W. Szyma{\' n}ski. Cuntz-Krieger 
algebras of infinite
graphs and matrices, preprint, Univ. of Newcastle, 1999.




\bibitem{rieffel:90} M.A.~ Rieffel.
\newblock Proper actions of groups on $C^*$-algebras,
\newblock {\em Mappings of Operator Algebras}, Birkh{\" a}user
Verlag, 1990, pages 141--182.

\bibitem{rieff99} M.A. Rieffel. Integrable and proper actions on 
$C^*$-algebras, and
square-integrable representations of groups, preprint, Univ. of 
California, Berkeley, 1999.

\bibitem{st} J.~Stillwell.
\newblock {\em Classical Topology and Combinatorial Group Theory,}
\newblock {Graduate Texts in
Math.}, vol. 72, Springer--Verlag, 1980.

\end{thebibliography}
\end{document}